\documentclass{article}\usepackage{amsmath}
\usepackage{latexsym}
\usepackage{amssymb}
\usepackage{amscd}

\newcommand{\ad}{\operatorname{ad}}

\newcommand{\EP}{Euler-Poincar\'e\ }
\newcommand{\h}{\hbar}
\newcommand{\cH}{\mathcal{H}}

\newcommand{\C}{\mathbb{C}}
\newcommand{\GM}{Groenewold-Moyal\ }

\newcommand{\Wqp}{W_{\text{qp}}}
\newtheorem{theorem}{Theorem}
\newtheorem{corollary}{Corollary}

\title{On the cohomology of the Weyl algebra, the quantum plane, and the q-Weyl algebra}
\author{M. GERSTENHABER and A. GIAQUINTO}

\begin{document}
\maketitle

\vspace{-7mm}
\date{}
{\noindent \textit{Department of Mathematics, University of
Pennsylvania, Philadelphia, PA 19104-6395 U.S.A., email:mgersten@math.upenn.edu} \linebreak[0]
U.S.A.\\ \noindent \textit{Department of Mathematics, Loyola University
Chicago, Chicago, IL 60626-5385 U.S.A., email:tonyg@math.luc.edu}}

\begin{abstract}\noindent Deformation theory can be used to compute the cohomology of a deformed algebra with coefficients in itself from that of the original. Using the invariance of the Euler-Poincar\'e characteristic under deformation, it is applied here to compute the cohomology of the Weyl algebra, the algebra of the quantum plane, and  the q-Weyl algebra.  The behavior of the cohomology when $q$ is a root of unity may encode some number theoretic information.
\end{abstract}

Algebraic deformation theory can be used to compute the cohomology of a deformed algebra with coefficients in itself when that of the original is known. This is illustrated here with the (first) Weyl algebra $W_1 = k\{x,y\}/(xy-yx-1)$, the algebra $W_{\text{qp}} = k\{x,y\}/(xy-qyx)$ of the quantum plane, and the q-Weyl algebra $W_q=k\{x,y\}/(xy-qyx-1)$,  where $k$ is a field of characteristic zero and $k\{x,y\}$ is the noncommutative polynomial ring in two variables over $k$. All three arise from deformations of the polynomial ring in two variables, $k[x,y]$.  An essential tool will be the invariance of the \EP characteristic of a possibly infinite complex under deformation. 

When a complex is deformed the cohomology classes of the deformed complex are those of cocycles which can be ``lifted" from ones of the original modulo those which lift to coboundaries,  \cite{G:DefIV}.  In particular, a deformation of an associative algebra $A$ induces a deformation of the Hochschild cochain complex $C^*(A,A)$ of $A$ with coefficients in itself.  Deformation of a complex preserves its \EP characteristic, provided that characteristic is well-defined, even if the dimensions of the cochain groups are infinite. This was in effect shown in \cite{G:DefIV} by considering obstructions and jump cocycles but will be reproved here more directly.

The Hochschild cohomology of an associative algebra $A$ with coefficients in itself, usually denoted $H^*(A,A)$,  will here generally be denoted simply $H^*(A)$ or just $H^*$, when $A$ is understood.  Sridharan \cite{Sridharan:Filtered} proved that $H^n(W_1) = 0, \,n\ge1$. (The result is sometimes attributed to the later extensive treatise of Dixmier \cite{Dixmier:Weyl}.) Although $W_1$ is not Koszul it has a `twisted' Koszul resolution from which its  cohomology and homology can be  computed. The latter  does not vanish identically cf e.g. \cite{Kassel:HomologyTrieste}. The cohomology of  $W_{\text{qp}}$ can be computed using Koszul resolutions, cf. \cite{Wambst:Hochschild}, \cite{Witherspoon:Hochschild}.  While $H^*(W_q)$ may be previously unknown, our main purpose  is not its computation but to illustrate some points of deformation theory.

The Weyl algebra $W_1$ can be exhibited in various ways as a
deformation of the polynomial ring $k[x,y]$, of which the
simplest is the {\em normal form}: Letting ``$*$'' denote the deformed
multiplication set
 $$a*b = ab
+\hbar\partial_xa\cdot\partial_yb +\frac{1}{ 2!}\hbar^2
\partial_x^2a\cdot\partial_y^2b + \dots =
m[\exp(\hbar\partial_x\otimes\partial_y)(a\otimes b)],$$ 
where on
the right $m$ denotes the original multiplication and $\hbar$ is a variable. That this defines an associative multiplication is a special case of the fact that if $\phi$ and $\psi$ are commuting derivations of an
algebra $A$ of characteristic zero, then setting $a*b =
m\exp(\hbar\phi\otimes\psi)(a\otimes b)$ defines a deformation of
$A$, \cite{G:DefIV}; the $2$-cocycle $\phi\smile\psi$ is its infinitesimal. In the special case where $\phi = \partial_x, \psi = \partial_y$ this is frequently attributed to Moyal
\cite{Moyal} but the idea is implicit earlier in Groenewold
\cite{Groenewold}; we will refer to these as  Groenewold-Moyal deformations. More generally, given commuting derivations
$\phi_1,\dots, \phi_r, \psi_1,\dots,\psi_r$ setting $a*b =
m\exp(\hbar\sum_i\phi\otimes\psi_i)(a\otimes b)$ defines a deformation
Also, since
$\phi\otimes\psi$ is cohomologous to $\phi \wedge \psi = (1/2)(\phi
\otimes \psi - \psi \otimes \phi)$ (the cup product is supercommutative in cohomology) using the latter as infinitesimal gives an isomorphic algebra.  That these formulas in fact produce associative deformations is a special case of \cite[Theorem 2.1]{GiaquintoZhang:Bialgebra}. Using either the normal or the preceding skew form for the star multiplication, and denoting by $[-,-]_*$  the commutator in this multiplication, one has $[x,y]_* = \hbar$. Taking $\hbar=1$ gives the Weyl algebra. In this case, specializing $\hbar$ to any non-zero element gives, up to isomorphism, the same algebra. This is characteristic of a \emph{jump deformation} (cf. \S\ref{jump}), that is, one which can be roughly described as not depending on the deformation parameter once that parameter is not zero.  (For jump  deformations, lifting of cocycles, and obstructions to their liftings, in the algebraic context, cf \cite{G:DefIV}.)  A deformation of an algebra $A$ generally has as underlying module $A[[h]]$ (power series over $A$) and its  coefficient ring is $k[[\hbar]]$ but jump deformations can frequently be defined  on $A[\hbar]$ with coefficient ring $k[\hbar]$.  The deformation parameter can then be specialized to an element of the ground field.

 When a deformation of an algebra $A$ is denoted by $A_{\hbar}$ it will be understood that this is the `generic' case, i.e., that $\hbar$ is a variable over the ground field.

\section{Deformation and cohomology} The Hochschild cochains of an algebra $A$ with coefficients in itself admit a bracket product \cite{G:DefIV}  (cf \cite[\S7]{GS:Monster}) in which the the product of cocycles is a cocycle and that of a cocycle with a coboundary is a coboundary; it is therefore defined on cohomology. If $F^i$ and $F^j$ are cochains of $A$ of dimensions $i$ and $j$, respectively, then that of $[F^i, F^j]$ is $i+j-1$. In low dimensions one has the following. The bracket of a one-cochain $F^1$ with a two-cochain $F^2$ is the two-cochain defined by $[F^1, F^2](a,b) = F^1(F^2(a,b))- (F^2(F^1(a),b) + F^2(a, F^1(b)))$. If $D, D_1, D_2$ are derivations  this gives $[D, D_1\smile D_2] = [D, D_1] \smile D_2 + D_1 \smile [D, D_2]$ where the brackets are now just the usual commutators. If $c\in A$  then $[F^2,c]$ is defined by $[F^2,c](a) = F^2(c,a)-F^2(a,c)$ and $[F^1,c]$ is just the element $F^1(c)$ of the coefficient ring $k$. Denoting the multiplication in $A$ by $m$, the Hochschild coboundary of a cochain $z$ is given by $\delta z = -[z,m]$. 

When $m$ is deformed to a star product $m_{\hbar} =m +\hbar m_1 +\hbar^2m_2 +\dots$ an $n$-cocycle $z$ is called {\em liftable} (cf \cite{G:DefIV}) to a cocycle of the deformed algebra $A_{\hbar}$ if there is a cocycle of $A_{\hbar}$  of the form $z_{\hbar} = z+\hbar z_1 +\hbar^2 z_2 +\dots$, where the $m_i$ are 
$k$-bilinear maps from $A\times A$ to $A$ which are tacitly extended to be bilinear over $k[[\hbar]]$ and the $z_i$ are $n$-cochains of $A$ similarly extended.  Since the criterion that $z_{\hbar}$ be a cocycle of $A_{\hbar}$ is that $[z_{\hbar}, m_{\hbar}] = 0$, expanding the preceding in powers of $\hbar$ shows that this is equivalent to the infinite sequence of conditions $[z, m_1] = \delta z_1,  [z, m_2] + [z_1,m_1] = \delta z_2, \dots$. Here $[z,m_1]$ is guaranteed to be a cocycle and the first of the sequence of equations asserts that its cohomology class is the primary obstruction to the existence of some $z_1$. If this class vanishes and a $z_1$ has been chosen then the second equation gives the obstruction to the existence of $z_2$, etc. . If the dimension of $z$ is $n$ then the obstructions lie in dimension $n+1$, so if the cohomology in that dimension vanishes then all obstructions vanish and $z$ is liftable.  The group $H^n(A_{\hbar})$ is isomorphic to the space of liftable $n$-cocycles of $A$ modulo the subspace of those that lift to coboundaries in $A_{\hbar}$ (but see the next paragraph). Coboundaries are always liftable and lift to coboundaries, so lifts are well-defined on cohomology classes.  Also, if two cocycles lift then their cup product can be lifted to the cup product of their respective lifts, and similarly for bracket products. If an algebra $A$ has $H^n(A) = 0$ for some $n$ then the same will be true for any deformation of $A$ as long as the deformation parameter is a variable, but this may fail when it is specialized to some value in the ground field. The dimensions of the cohomology groups need not be the same at different values of the deformation parameter. 

More precisely, a cocycle $F$ which is not a coboundary of an algebra $A$ can lift to a coboundary of the deformed algebra $A_{\hbar}$ only after extension of the coefficients to include $\hbar^{-1}$. For if $\delta_{\hbar}$ is the deformed coboundary operator then we can not have a cochain $f_{\hbar} = f + \hbar f_1 + \dots$ with $\delta_{\hbar}f_{\hbar} = F_{\hbar}$ where $F_{\hbar}$ is the lift of $F$ since 
setting $\hbar = 0$ would show that $F$ is already a coboundary in $A$.  If $F$ lifts to a coboundary then one must have $\delta f= 0$, and possibly also $\delta f_i = 0, i = 1, \dots, r-1$ for some $r$, and $\delta_{\hbar}f_{\hbar} = h^rF_{\hbar}$, so $F_h$ is not a coboundary but $h^rF_{\hbar}$ is; when considered over $k[[\hbar]]$ the cohomology of $A_{\hbar}$ then has $\hbar$ torsion. This technical difficulty disappears with a jump deformation when one specializes $\h$ to 1.

 \begin{theorem}\label{central}  A derivation which lifts to a coboundary annihilates any central element which can be lifted to a central element. 
\end{theorem}
\noindent\textsc{Proof}. Let $A_{\hbar}$ be a deformation of an algebra $A$ and
suppose that a derivation $D$ of $A$ lifts to a derivation $D_{\hbar}$ of $A_{\hbar}$ which is inner in $A_{\hbar}$, i.e., that there is a $\xi_{\hbar} = \xi + \hbar \xi_1 + \dots$ such that $h^iD_{\hbar} = \operatorname{ad}\xi_h$ for some $i$.  Suppose also that a central element  $c$ of $A$ lifts to a central element $c_{\hbar}$ of $A_{\hbar}$.  Then $h^iD_{\hbar}c_{\hbar} = 0$ identically, so $D_{\hbar}c_{\hbar} = 0$, whence $Dc = 0$. $\Box$

\medskip
The theorem will be used to show that certain lifts of derivations in fact do not become inner.

\section{Invariance of the Euler-Poincar\'e characteristic}\label{EP} 
A finite cochain complex of finite dimensional vector spaces 
\begin{equation} \label{complex}
\begin{CD}
V^{\bullet} :\quad 0\to  V^0 @>\delta>> V^1 @>\delta >>\cdots @>\delta >> V^{n-1} @>\delta >> V^n \to 0
\end{CD}
\end{equation}
in principle has two \EP characteristics, its dimensional one, $\chi_d = \sum (-1)^i \dim V^i$, and its homological one, $\chi_d = \sum (-1)^i \dim H^i$, where $H^i$ is the $i$th cohomology group of $V^{\bullet}$.  Euler (1706-1783) observed  that for polyhedra these coincide. That they always do is sometimes called the {\em \EP principle}, which being needed below, we reprove. 
\begin{theorem} With the above notation,
$$\chi_d = \chi_h.$$
\end{theorem}
\noindent\textsc{Proof}.
Recall that over any coefficient ring a short exact sequence of cochain complexes 
$$0\to C^{\bullet} \to V^{\bullet} \to B^{\bullet} \to 0$$
 gives rise to a long exact sequence in cohomology,
\begin{equation}\label {long}
\cdots \to H^i(C) \to H^i(V) \to H^i(B) \to H^{i+1}(C) \to \cdots.\end{equation}
Since by hypothesis the complex terminates so does the sequence (\ref{long}). Now apply alternating signs to the dimensions of the groups in (\ref{long}). Since the dimensional \EP characteristic of a finite exact complex clearly vanishes, (\ref{long}) then says precisely that 
\begin{equation}\label{chi}\chi_h(C) -\chi_h(V) +\chi_h(B) = 0.\end{equation}
This holds whatever the coboundary operators $\delta$ may be; setting them all identically equal to zero shows also (the obvious fact) that the same holds for the dimensional \EP characteristics.
If, further,  $C^{\bullet}$ is exact, then  its homological and dimensional  \EP characteristics both vanish so (\ref{chi}) implies that the characteristics, both homological and dimensional, of $V^{\bullet}$ and $B^{\bullet}$ coincide.  Now let $V^{\bullet}$ denote our original complex (\ref{complex})  and construct a subcomplex $C^{\bullet}$ of  $V^{\bullet}$ by choosing, for each $i$, an arbitrary vector space complement $W^i$ to the space of $i$-cocycles and setting $C^i = W^i + \delta V^{i-1}$. This subcomplex is exact and the quotient complex, $B^{\bullet} =V^{\bullet}/C^{\bullet}$ consists, up to isomorphism,  of just the cohomology groups of $V^{\bullet}$ with zero coboundaries. The dimensional and homological \EP characteristics of this quotient are formally identical, but since they coincide with those of $V^{\bullet}$, so are those of $V^{\bullet}$, as asserted.$\Box$

\medskip
The following, while not needed here, is a classic consequence.
 \begin{corollary}  The \EP characteristic of a finite dimensional Lie algebra $L$ with coefficients in a finite dimensional module always vanishes. 
\end{corollary}
\noindent\textsc{Proof}. The Chevalley-Eilenberg cochains are alternating, so if $\dim{L} = n $ then the dimension of the $i$th cochain group is $\binom{n}{i}$ times the dimension of the module, but the alternating sum of the binomial coefficients vanishes.  $\Box$
\smallskip

This is independent of the characteristic. The first proof  in full generality seems to be due to S. Goldberg \cite{Goldberg:Lie}; for a discussion of its history and an extension to the case of Lie superalgebras see Zusmanovich \cite{Zusmanovich:superLie}. Note that likewise the \EP characteristic of a non-trivial compact connected Lie group must vanish, else left multiplication by any group element would always have a fixed point by Lefschetz' theorem.

 It follows from the \EP principle that even if the boundary operators in a finite cochain complex are changed, the Euler-Poincar\'e characteristic remains constant. This  can not extend in full generality to an  infinite complex of infinite dimensional vector spaces.  Invariance of the \EP characteristic does hold, however for a deformation of a complex $V^{\bullet}$ which is not necessarily finite but which does have a well-defined \EP characteristic, and the proof is almost identical.
Label the coboundary operators more precisely as $\delta^{(n)}:V^n\to V^{n+1}$ so that we have $\delta^{(n+1)}\delta^{(n)}=0$.  A deformation of the complex is then a sequence of formal operators $\delta^{(n)}_{\hbar}:V^n[[\hbar]]\to V^{n+1}[[\hbar]]$ where $\delta^{(n)}_{\hbar}$ has the form
$\delta^{(n)}_{\hbar} = \delta ^{(n)} +\hbar\delta^{(n)}_1 +\hbar^2\delta^{(n)}_2 + \cdots$, each $\delta^{(n)}_i$ being a linear operator $V^n \to V^{n+1}$.  However, for this to be meaningful one must first replace each vector space $V^n$ by $V^n[[\hbar]]$ (power series in $\hbar$ with coefficients in $V^n$). This is a module over $k[[\hbar]]$, where $k$ was the original coefficient ring, which for the present is assumed to be a field.  With the obvious extension of the coboundary operators to be linear over $k[[\hbar]]$ the extension of the original complex can be denoted $V^{\bullet}[[\hbar]]$ and it is this which is being deformed. When necessary to continue working over a field one can adjoin $\hbar^{-1}$, the resulting complex $V^{\bullet}((\hbar))$ then has the field of formal Laurent series $k((\hbar))$ as coefficients.  Note that when $V$ is a finite dimensional vector space over a field $k$ (or more generally, a free module of finite rank)  then $V[[\hbar]]$ is identical, as a module over $k[[\hbar]]$ with $V\otimes_k k[[\hbar]]$ but when the dimension (or rank) of $V$ is infinite then the former is much larger. However, viewing $V[[\hbar]]$ and $V((\hbar))$ as having a uniform complete topology with the $\hbar^iV[[\hbar]]$ and $h^iV((\hbar))$, respectively, as neighborhoods of 0 (where in the latter case, $i$ may be negative), then any original basis of $V$ will continue to serve as a basis. If one interprets $\otimes$ as a suitably completed tensor product then one still has $V[[\hbar]] = V\otimes_k k[[\hbar]]$. 

\begin{theorem} \label{invariance} Suppose that a complex $V^{\bullet}$ of vector spaces over a field $k$ is not necessarily finite but has cohomology groups $H^n(V)$ which are all finite dimensional and vanish for sufficiently large $n$. Then its \EP characteristic remains well defined and is invariant under a deformation $V^{\bullet}_{\hbar}$ of the complex.  One then has $\dim H^n(V^{\bullet}_{\hbar}) \le \dim H^n(V^{\bullet}$), all $n$. 
\end{theorem} 
\noindent\textsc{Proof}.
Choose in each $V^n$ a subspace $W^n$ which is a      
vector space complement to the cycles $Z^n$ of $V^n$. Then no element of $W^n((\hbar))$ can be annihilated by $\delta^{(n)}_{\hbar}$, for multiplying by a suitable power of $\hbar$ we may assume, without loss of generality, that such an element has the form $v_0 + \hbar v_1 +\hbar^2v_2 + \cdots$  with $v_i \in W, v_0 \ne 0$; setting $\hbar = 0$ then gives $\delta^{(n)} v_0 = 0$, a contradiction.   

Set $C^n = W^n + \delta^{(n-1)}W^{n-1}$ and $C^n_{\hbar} = W^n((\hbar)) + \delta^{(n-1)}_{\hbar}W^{n-1}((\hbar))$.  From the foregoing, the sum is direct in both cases, the $C^n$ form an exact subcomplex $C^{\bullet}$ of $V^{\bullet}$ with the original coboundary operators and the
$C^n_{\hbar}$ likewise form an exact subcomplex of the deformed complex  $V^{\bullet}_{\hbar}$. Since, $V^{\bullet}$ had a well-defined \EP characteristic, there is for every $n$ a finite dimensional subspace $\cH^n$ whose dimension is that of $H^n(V^{\bullet})$ (and 
hence vanishes for sufficiently large $n$)  such that $V^n = W^n + \delta^{(n-1)}W^{n-1} +\cH^n$, the sum being direct. It follows that the sum $W^n((\hbar)) + \delta^{(n-1)}_{\hbar}W^{n-1}((\hbar)) + \cH^n((\hbar))$ is also direct since any relation multiplied by a suitable power of $\hbar$ will give one with 
no negative powers of $\hbar$ and not divisible by $\hbar$, from which one would get a relation in $W^n + \delta^{(n-1)}W^{n-1} + \cH^n$ by setting $\hbar = 0$. Moreover, one must have $W^n((\hbar)) + \delta^{(n-1)}_{\hbar}W^{n-1}((\hbar)) +\c H^n((\hbar)) = V^n((\hbar))$  since completeness implies that any set of elements of $V^n((\hbar))$ whose leading coefficients form a basis of $V$ must itself form a basis for $V^n((\hbar))$ and taking leading coefficients of the elements of  $W^n((\hbar)) + \delta^{(n-1)}_{\hbar}W^{n-1}((\hbar)) +\cH^n((\hbar))$ just gives $W^n + \delta^{(n-1)}W^{n-1} +\cH^n = V^n$.  

With the above definitions, in (\ref{long}) now replace $C^{\bullet}$ by $C^{\bullet}_{\hbar}$, $V^{\bullet}$ by $V^{\bullet}_{\hbar}$ and $B^{\bullet}$ by $V^{\bullet}_{\hbar}/C^{\bullet}_{\hbar}$.  Since $C^{\bullet}_{\hbar}$ is exact we have $H^i(V^{\bullet}_{\hbar}) \cong H^i(V^{\bullet}_{\hbar}/C^{\bullet}_{\hbar})$, so the \EP characteristic of $V^{\bullet}_{\hbar}$  is the same as that of the finite complex $V^{\bullet}_{\hbar}/C^{\bullet}_{\hbar}$ in which the $n$th vector space is isomorphic to $\cH^n$. The homological \EP characteristic of the quotient is the same as its dimensional \EP characteristic, but the latter is just the \EP of $V^{\bullet}$, proving the assertion that the \EP of the deformed complex is the same as that of the original. 

Before deformation, the coboundary operator in the quotient complex $V^{\bullet}/C^{\bullet}$ was identically zero. That need no longer be the case after deformation, so the $n$th cohomology group of $V^{\bullet}_{\hbar}$ can not have dimension greater than that of $\cH^n$;  it follows that the dimensions of the cohomology groups can not increase under deformation, as asserted. $\Box$
\medskip

Deformation of an associative algebra induces a deformation of its Hochschild complex, so the theorem applies to show, in particular, that the \EP characteristic, if well-defined, remains unchanged. This was previously known by examining obstructions to liftings of cocycles. 

\section{Gradations and group actions}\label{gradations}
When a complex of vector spaces is a direct sum of subcomplexes finite or not, each of which has a well-defined \EP characteristic, then a deformation which respects the direct sum decomposition will leave the \EP characteristic of each subcomplex unchanged. For example, the cohomology of a commutative algebra of characteristic zero has a Hodge-type decomposition introduced in \cite{GS:Hodge}. The Euler-Poincar\'e characteristic of each component may be well-defined even when the whole is not, and is preserved by deformation. 

If an algebra $A$ has a non-negative integral gradation $A=\bigoplus_{i=0}^{\infty}A^i$, where $A^iA^j \subseteq  A^{i+j}$ then  define an $n$-tuple $(a_1,\dots,a_n)$ of homogeneous elements to have degree the sum of the degrees of its components, and define a Hochschild $n$-cochain $f\in C^n(A)$ to be homogeneous of degree $r$ if $f(a_1,\dots,a_n)$ is homogeneous of degree $\sum\deg a_i+r$ whenever the $a_i$ are themselves homogeneous. Each homogeneous part of the cochain complex is then carried into itself and so is a subcomplex of the total Hochschild complex.  Suppose now that one has a deformation of $A$ which respects the grading. This will induce a deformation of each of the homogeneous subcomplexes, so the \EP characteristic of each, if defined, will be preserved.  If almost all of the subcomplexes composing the direct sum have vanishing Euler-Poincar\'e characteristic, then we can assign to the total complex a characteristic equal to the sum of those of the subcomplexes and this will be constant under any deformation respecting the decomposition.  

It is generally not meaningful to specialize the parameter $\hbar$ in the definition of the deformation of a complex $V^{\bullet}$ to an element of $k$ unless either the vector spaces involved carry a topology in which convergence of the series defining the deformation is meaningful for a certain range of the parameter or $\delta_{\hbar}v$ is a polynomial in $\hbar$ for every vector $v$. While $V^{\bullet}$, and hence also $V^{\bullet}_{\hbar}$, may have had well-defined (and by the above, equal) \EP characteristics, it does not follow that this holds after specialization: At a particular value $\hbar_0$  of $\hbar$ (assuming that specialization of $\hbar$ to $\hbar_0$ is meaningful) some cohomology groups may become infinite dimensional and there may be infinitely many non-zero groups. In the other direction, it can happen that $H^i(V^{\bullet}_{\hbar}) \ne 0$ for some $i$ while  $H^i(V^{\bullet}_{\hbar_0}) = 0.$  For example, the $q$-Weyl algebra $W_q$ has both $H^1(W_q) \ne 0$ and $H^2(W_q) \ne 0$ whenever $q \ne 1$ but these groups vanish when $q=1$; the Weyl algebra  $W_1$ possesses no non-trivial cohomology.  This kind of anomalous behavior can occur only for infinite dimensional algebras, but we are dealing with such here.

Despite the foregoing,  if the generic deformed complex $V^{\bullet}_{\hbar}$ can be viewed also as a deformation of the specialized deformed complex $V^{\bullet}_{\hbar_0}$  then the \EP characteristic of the specialized deformed complex certainly must remain well-defined and equal to that of the original complex.  Moreover, in that case if $H^i(V^{\bullet}_{\hbar}) \ne 0$ for some $i$ then also $H^i(V^{\bullet}_{\hbar_0}) \ne 0$.  This is so, for example, when we have deformations of $V^{\bullet}$ which are naturally parameterized by a group, as is the case with Groenewold-Moyal deformations.  While we do not use its full generality here, recall from 
 (\cite[Theorem 2.1]{GiaquintoZhang:Bialgebra}) the following:  Suppose that the coefficient ring $k$ contains $\mathbb{Q}$, that $B$ is a commutative bialgebra, and that $P$ is its space of primitive elements. If $r$ is any element of $P\otimes P$ then 
$$\exp(\hbar r) =\sum_{i=0}^{\infty}\frac{\hbar^i}{i!}r^i = 1\otimes 1 + \hbar r + \frac{\hbar^2}{2!}r^2 + \cdots$$
is a `universal deformation formula'.  Observe that there is an action of the additive group, since $\exp(\hbar_1 r)\exp(\hbar_2 r) = \exp((\hbar_1 +\hbar_2)r)$.  But this says that if we have applied such a deformation and specialized $\hbar$ to some $\hbar_0$ then we can return to the generic value of $\hbar$ by applying $\exp((\hbar-\hbar_0)r)$.  
Specializing to $\hbar = 0$ in the foregoing returns the original algebra, which may be paradoxical since one may start with a polynomial ring and deform it to a non-commutative algebra. But note that while a generic deformation of a non-commutative algebra must remain non-commutative, a specialization of that deformation may well become commutative and even return us to the original algebra. This leads to the question of what is the modular group of a deformation problem, cf \cite[p.91]{G:DefI}. If we have a Groenewold-Moyal deformation, or more generally, a deformation parameterized by a group action, then as above the generic deformation is also a deformation of any specialization. Therefore, if the original algebra and the specialization have well-defined \EP characteristics then they must coincide since the generic deformation is a deformation of both.  The condition that $A_{\hbar_0}$ be well-defined for all $\hbar_0$, while exceedingly strong, is satisfied for most jump deformations (next section).

The polynomial ring $k[x,y]$ is bigraded; the bidegree of $x^ry^s$ is $(r,s)$.  Since $H^0 $ is the center it is all of $k[x,y]$; its part of degree $(r,s)$ consists of the multiples of $x^ry^s$.  Again, because $k[x,y]$ is commutative, there are no inner derivations so $H^1$ consists of all derivations. These are likewise bigraded, the degree of $x^ry^s\partial_x$ being $(r-1,s)$ and that of $x^ry^s\partial_y$ being $(r,s-1)$. By the theorem of Hochschild-Kostant-Rosenberg \cite{HKR}, $H^2$ is isomorphic to 
$H^1\wedge H^1$ which is spanned by elements of the form $D_1\wedge D_2 = (1/2)(D_1\smile D_2 - D_2\smile D_1)$ where $D_1, D_2$ are derivations of $k[x,y]$, and $H^i =0$ for $i\ge 3$. Note that in a commutative algebra $A$ of characteristic $\ne 2$ every $2$-cocycle can be written as a sum of its symmetric and skew parts, each of which is again a cocycle. The coboundary of every $1$-cochain is always symmetric. In general a symmetric $2$-cocycle need not be a coboundary but no skew $2$-cocycle can be cohomologous to zero, so each skew $2$-cocycle is the unique representative of its cohomology class. Thus $H^2(k[x,y])$ is spanned by the classes of the $2$-cocycles $x^ry^s\partial_x\wedge\partial_y$, and these are linearly independent. Each has a bidegree, namely 
$(r-1,s-1)$; symmetric $2$-cocycles must be coboundaries.  In any deformation there can be no obstruction to lifting an element of $H^2$ since $H^3$ = 0.  The Euler-Poincar\'e characteristic $\chi$ of $k[x,y]$, which is not initially well-defined, decomposes as follows. Denote by $H^n_{r,s}$ the part of $H^n$ of bidegree $(r,s)$, by $h^n_{r,s}$ its dimension, and set $\chi_{r,s} = \sum (-1)^i h^n_{r,s}$. Note that $h^n_{r,s} = 0$ for $n\ge 3$.   For every bidegree $(r,s)$ with $r,s \ge 0$ one has $h^0_{r,s} = 1, h^1_{r,s} = 2, h^2_{r,s} = 1$ so $\chi_{r,s} = 0$. For $s\ge 0$, $h^0_{-1,s} = 0,  H^1_{-1,s}$ is spanned by the class of $y^s\partial_x$  and $H^2_{-1,s}$ is spanned by the class of $\partial_x\wedge y^s\partial_y$ so $\chi_{-1,s} = 0$, and similarly $\chi_{r, -1}= 0$ for $r \ge 0$. The lowest bidegree that can occur is $(-1,-1)$, which happens only for the class of $\partial_x\wedge\partial_y$, so $\chi_{-1,-1} = 1$. Thus, although the Euler-Poincar\'e characteristic of $k[x,y]$ is initially undefined since the dimensions of the groups are infinite, if we define it to be the sum of those of its bigraded parts then the Euler-Poincar\`e characteristic becomes precisely $1$. 

For more detailed treatments of the deformation of complexes, including analysis of obstructions, cf. \cite{GerstWilker:DefV}, \cite{ArmstrongUmble:Obstructions}.

\section{Jump deformations}\label{jump} Intuitively, a jump deformation of an algebra $A$  is one which does not depend on the deformation parameter once that parameter is different from zero. The precise definition in the algebraic context, introduced in \cite{G:DefIV}, is this. Suppose that the deformed multiplication on $A_{\hbar}$ (which is now defined over $k[[\hbar]]$) is given by
$$a\star_{\hbar} b = ab +\hbar m_1(a,b) + \hbar^2m_2(a,b) + \dots.$$
Then  we can  further deform it by ``stretching'' the deformation parameter: Introduce another parameter, $t$ and replace $\hbar$ by $(1+t)\hbar$. This gives an algebra $A_{(1+t)\hbar}$ with multiplication
$$a\star_{(1+t)\hbar}b =  ab +(1+t)\hbar m_1(a,b) + (1+t)^2\hbar^2m_2(a,b) + \dots,$$
which is now defined over $k[[\hbar]][[t]]$.
A {\em jump deformation} is one such that $A_{(1+t)\hbar}$ is a trivial deformation of $A_{\hbar}$; a linear stretching of the deformation parameter does not change the deformed algebra.  This may seem weaker than the intuitive definition, but it is enough. The infinitesimal of this second trivial deformation, however, is just the original $m_1$, which implies that $m_1$ has become a coboundary in the deformed algebra:
\begin{theorem}
 The infinitesimal of a jump deformation lifts to a coboundary in the deformed algebra.
$\Box$
\end{theorem}
 A jump deformation thus  ``annihilates its own infinitesimal.'' The converse need not be true since this is only a first order property but raises the following question: Suppose that a deformation has the property that its infinitesimal lifts to a coboundary; does there then exist a jump deformation with the given infinitesimal?  

If we have a jump deformation  and extend coefficients to include $\hbar^{-1}$ then setting $t= \hbar^{-1}-1$ shows that $A_{\hbar}$ is isomorphic to $A_1$, whose multiplication, providing the following formula is meaningful, is given by
$$a\star_1 b = ab +m_1(a,b) +m_2(ab) + \dots.$$
This is so in some of the most important examples, namely those where there is a topology in which the series converges and those, including those here, where the series terminates for any fixed $a$ and $b$.  In the latter case the deformed algebra could already have been defined over $k[\hbar]$ independent of whether the deformation was a jump. But in any case where $A_1$ is defined it is clear that the deformed algebra  does not depend on the deformation parameter once the latter is different from zero 

\section{The Weyl algebra} This section presents a brief proof of Sridharan's theorem on the vanishing of the cohomology of the Weyl algebras.  The algebra $W_{\h} = k[\h]/(xy-yx-\h)$  is readily seen to be a jump deformation of $k[x,y]$.  One has $W_{\h} \cong W_1 \otimes_k k[\h]$, so one often says simply that $W_1$ is a jump deformation of $k[x,y]$. The $i$th Weyl algebra is the $i$th tensor power of $W_1$.

\begin{theorem}  If $W$ is the $i$th  Weyl algebra, then $H^0(W) = k$ and $H^n(W) = 0$ for $n \ge 1$.
\end{theorem}
\noindent\textsc{Proof}.
It is sufficient to prove this for the first Weyl algebra $W_1$. The first assertion is simply that the center is reduced to $k$, which is evident. Since tensoring with $k[\hbar]$ will not destroy any cohomology, it is sufficient to prove the remaining assertions for 
$W_{\hbar} = k\{x,y\}/(xy-yx-\hbar)$, which may be obtained by deformation of $k[x,y]$ using the normal form of the $*$ product
$$a*b = ab +\hbar\partial_xa\cdot \partial_yb + (t^2/2!)(\partial_x)^2a\cdot (\partial_y)^2b + \dots.$$ The infinitesimal $m_1$ of this deformation is $\partial_x\smile \partial_y$.  Since $H^i(k[x,y]) = 0$ for $n >2$ the same must be true by Theorem \ref{invariance} for $W_{\hbar}$, so Sridharan's theorem reduces to showing that $H^i(W_{\hbar}) = 0$ also for $i = 1,2$. 

Every derivation of $k[x,y]$ is of the form $D = a(x,y)\partial_x + b(x,y)\partial_y; \, a, b\in k[x,y]$. The primary obstruction to lifting $D$ is 
\begin{multline*}
 [D, \partial_x\smile \partial_y] = [D, \partial_x]\smile\partial_y +\partial_x \smile [D, \partial_y] =\\
-(a_x\partial_x\smile\partial_y + b_x\partial_y\smile\partial_y)
  -(a_y\partial_x\smile\partial_x + b_y\partial_x\smile\partial_y)
\end{multline*}
  where we have written $a_x$ for $\partial_x a$, etc.. Now $ \partial_x\smile\partial_x = -(1/2)\delta(\partial_x)^2$, and so is a coboundary, as is $ \partial_y\smile\partial_y$, but $ \partial_x\smile\partial_y$ is not a coboundary, so the primary obstruction to lifting $D$ vanishes if and only $a_x = - b_y$. This is equivalent 
to the existence of $c\in k[x,y]$ with $c_x = a, c_y = -b$. But then $(1/\hbar) \operatorname{ad} c$ is a lift of $D$ to $W_{\hbar}$, for if $w\in k[x,y]$ then the series for $[c,w]_*$ (commutator in the deformed product) begins with 
$\hbar(c_xw_y - c_yw_x)$.  As every derivation of $k[x,y]$ which can be lifted must have vanishing first obstruction and every such lifts to a coboundary one has $H^1(W_{\hbar})= 0$. 

For the second cohomology, every two-cocycle of $k[x,y]$ (which is always liftable since its obstruction would be a three-cocycle, hence a coboundary)  is cohomologous to exactly one of the form $a(x,y)\partial_x\smile \partial_y$. These are linear combinations of ones of the form $x^ry^s\partial_x\smile \partial_y$ so it is sufficient to show that the latter  lift to coboundaries.  But from the foregoing $y^s\partial_x$ and $x^r\partial_y$ do lift, each to a coboundary, so their cup product $x^ry^s\partial_x\smile \partial_y$ lifts to the cup product of these respective coboundaries which is again a coboundary, hence $H^2(W_1) = 0$ as well. $\Box$
\medskip

The Euler-Poincar\'e characteristic (which was not well-defined initially but defined as in \S2) has been preserved by the deformation, for that of the Weyl algebra (which is well-defined from the start) is just the dimension of $H^0$, its center, i.e., $1$. In view of this example, one might try to define the \EP characteristic of an algebra with infinite dimensional or infinitely many non-zero cohomology groups as that of some generic deformation whose \EP characteristic is well-defined.  This actually assigns an \EP characteristic to a family of algebras but the problem is that an algebra may lie at the intersection of two families with different \EP characteristics. 

\section{The quantum plane}  The Groenewold-Moyal deformation of $k[x,y]$ with infinitesimal $x\partial_x\smile y\partial_y$ (or equivalently $x\partial_x\wedge y\partial_y$) is readily seen to give a star product where formally $x\ast y = e^{\h}y\ast x$. When $k = \C$ it is meaningful to specialize $\h$ to any complex value, so setting $q=e^{\hbar}$, we may view $W_{\text{qp}}$ with $q$ any non-zero complex number as a deformation of $k[x,y]$. It is not a jump deformation, for it depends essentially on the parameter $q$, but it is Koszul and has a projective resolution of length 2 (cf, e.g.  Kraehmer \cite{Kraehmer:Koszul}). It follows that $H^n =0$ for $n>2$ not only for generic $q$, where it follows from deformation theory, but for arbitrary $q\in \C$. (To show that $H^n(\Wqp ,M)$ vanishes for $n >2$ with $M$ any $\Wqp$ bimodule it is sufficient, and relatively simple, to show that $k$, considered as a trivial $\Wqp$ bimodule, has a projective resolution of length 2; again cf. \cite{Kraehmer:Koszul}.) Note that $W_{\text{qp}}$ preserves the bigrading of $k[x,y]$, so in the following we only have to consider homogeneous cocycles; also, the cohomology of an algebra is always a module over its center.

\begin{theorem}\label{qp} 1. When $q$ is a primitive $N$th root of unity then the center $H^0$ of $W_{\text{qp}}$ is generated by $x^N, y^N$ and is reduced to $k$ otherwise. 2. $H^1$ is a free module of rank 2 over the center with generators the classes  of $x\partial_x, y\partial_y$, which remain well-defined in $W_{\text{qp}}$; it is just the two dimensional vector space over $k$ spanned by these two classes if $q$ is not a root of unity. 3. $H^2$ is the sum of the $k$ subspace spanned by the class of a cocycle $z_{\text{qp}}$ which is a lift of the cocycle $\partial_x\wedge\partial_y$ of $k[x,y]$ and the free module over the center generated by the class of $x\partial_x\wedge y\partial_y$.  
\end{theorem}
\noindent\textsc{Proof.}  The first assertion is evident. As for the second,
the primary obstruction to lifting an arbitrary derivation $a\partial_x + b\partial_y, \, a,b \in k[x,y]$ to a derivation of $W_{\text{qp}}$ is the cohomology class of 
$$[a\partial_x + b\partial_y,  x\partial_x \smile y\partial_y] =
 [a\partial_x + b\partial_y, x\partial_x]\smile y\partial_y +
x\partial_x \smile[a\partial_x + b\partial_y,  y\partial_y].$$
Since $\partial_x\smile\partial_x$ and $\partial_y\smile\partial_y$ are coboundaries but $\partial_x\smile\partial_y$ is not, this is a coboundary if and only if the coefficient of $\partial_x\smile\partial_y$ vanishes.
Writing $\partial_xa = a_x$ and similarly for $b_y$, that coefficient is $ay - xa_xy + bx - yb_yx$. This vanishes if and only if  $(a/xy)_x = -(b/xy)_y$, which is precisely the condition that there exist a function $f(x,y)$ such that $a = xyf_y, b = -xyf_x$. Since $a$ and $b$ must be polynomials in $x$ and $y$, the function $f$ must itself be a polynomial plus some constant multiples of $\ln x$ and $\ln y$.  The latter two give, respectively, the derivations $-y\partial_y$ and $x\partial_x$, which simply lift to themselves since they commute with the infinitesimal of the deformation;  any such derivation of the original algebra  simply lifts to itself.  This can not be a coboundary unless it originally was one, since it is not a function of the deformation parameter. Therefore 
$x\partial_x$ and $y\partial_y$ are not inner derivations of the deformed algebra, something also easily seen directly.  All other homogeneous derivations of $W_{\text{qp}}$ have non-negative bidegrees and, up to constant multiple, are obtained by setting $f = x^ry^s$ with at least one of $r$ and $s$ strictly positive; this yields the derivation $D_{rs} = sx^{r+1}y^s\partial_x -rx^ry^{s+1}\partial_y$ of bidegree $r,s$. It follows that there is at most one derivation of $W_{\text{qp}}$ of bidegree $r,s$.  However, $\ad x^ry^s$ is also a homogeneous derivation of bidegree $r,s$, so $D_{rs}$ lifts to a coboundary unless  $\ad x^ry^s=0$, i.e., unless $x^ry^s$ is central. On the other hand, $H^1$, which contains $x\partial_x$ and $y\partial_y$ is always a module over the center and no multiple of any of these by a central element can be a coboundary by the foregoing. This yields the second assertion.

For the third assertion, note that since the bigrading is preserved so is the \EP characteristic $\chi_{r,s}$ of each homogenous part $H^*_{r,s}$ of bidegree $(r,s)$ of the cohomology.  As observed in \S\ref{EP} 
this is zero except for $(r,s) = (-1,-1)$. Further, $k[x,y]$ has no cohomology in any lower bidegree so the same must be true of  $H^2(W_{\text{qp}})$. Consider first the case of non-negative bidegree. If $H^0_{rs} = 0$, i.e., if $x^ry^s$ is not central, then from the preceding preceding,  $H^1_{r,s} = 0$ as well. It follows 
then that also $H^2_{r,s} = 0$ since there is no higher cohomology.  If, however, $x^ry^s$ is central, so $\dim H^0_{r,s} = 1$, then from the preceding $\dim H^1_{r,s} = 2$ so we must have $\dim H^2_{r,s} = 1$ in order that $\chi_{r,s}= 0$. The 2-cocycle $x\partial_x\wedge y\partial_y$ of $k[x,y]$ lifts to itself (being a product of 1-cocyles with the same property) and can not become a coboundary since it does not depend on the parameter $q$. It follows that the sum of all $H^2_{i,j}$ with $i,j\ge0$ is a free module of rank one over the 
center of $W_{\text{qp}}$ generated by the class of $x\partial_x\wedge y\partial_y$.  The same arguments show that $H_{-1,0} = H_{0,-1} = 0$. It remains only to consider $H_{-1,-1}$. The cohomology of $k[x,y]$ does have a two dimensional class of bidegree $(-1,-1)$, namely that of $\partial_x\wedge \partial_y$. This must lift to some non-trivial cocycle, which we denote by $z_{\text{qp}}$, in order that $\chi_{-1,-1}$ remain equal to 1,  proving the last assertion. $\Box$

\medskip 
The Koszul resolution in  \cite{Witherspoon:Hochschild} exhibits a class, there denoted $1\otimes x_1^*\wedge x_2^*$, corresponding here to the lift $z_{\text{qp}}$ of $\partial_x\wedge\partial_y$, with the unusual property (as can be seen from the foregoing) that all $x^{rN}y^{sN}z_{\text{qp}}$ with either $r$ or $s$ strictly positive must be coboundaries.  The underlying reason is that while $z_{\text{qp}}$ is the infinitesimal of the deformation of the quantum plane to the generic $q$-Weyl algebra  $W_q(\h) = k\{x,y\}/(xy-qyx-\h)$ (\S\ref{q-Weyl}), replacing $\h$ in the defining equation by, say, $x^n\h$ for any $n > 0$ reproduces the quantum plane: just replace $y$ by $y+x^{n-1}\h/(1-q)$.

\section{The $q$-Weyl algebra}\label{q-Weyl} 
The remarkable thing about the cohomology of the q-Weyl algebra is that it is practically identical with that of the quantum plane, both when $q$ is a root of unity and when it is not.  In all cases only two cohomology classes are lost. 
 
The $q$-Weyl algebra with deformation parameter a variable $\hbar$, given by $W_q(\hbar) =k[h]\{x,y\}/(xy-qyx-\hbar)$, is a jump deformation of $W_{\text{qp}}$; the usual $q$-Weyl algebra is $W_q(1)$. One has $W_q(\hbar) \cong W_q \otimes k[\h]$. The cohomology of $W_1$ vanishes in positive dimensions by Sridharan's theorem, so $W_1$ is absolutely rigid.  Thus $W_q$, whose cohomology is non-trivial, can not be a deformation of $W_1$ in the classical sense. The classic deformation theory of a \emph{single} algebra can not detect the change from $W_1$ to $W_q$ but the deformation theory of a diagram (presheaf) of algebras does; cf \cite{GG:DefsRigidAlgs}. This is a phenomenon of infinite dimensional algebras, where to understand the deformation of a single algebra one must consider also diagrams of algebras built from the original.  
Nevertheless, $W_q$ is a deformation of $W_{\text{qp}}$, from which we can determine its cohomology. The infinitesimal of this deformation is just the cocycle $z_{\text{qp}}$ of the previous section. Since this is a jump deformation it lifts to a coboundary; $z_{\text{qp}}$ is one of the two classes which is lost. As we continue to have $H^n(W_q)$ = 0 for $n>2$,  it is only necessary to compute the center $H^0$, $H^1$, and $H^2$ of $W_q$.

In the following we will denote $1+q+q^2+\cdots +q^{n-1}$ by $n_q$. For convenience we may write $r= 1/q$ and then $n_r = 1+r+\cdots +r^{n-1}$. By induction one then has 
\begin{equation}\label{induction}
yx^n - r^nx^ny = n_r(\h r)x^{n-1}.
\end{equation}

\begin{theorem}  If $q$ is a primitive $N$th root of unity then the center of $W_q$ is $k[x^N,y^N]$ and is reduced to $k$ otherwise. It is identical with the center of $W_{\text{qp}}$.
\end{theorem}
\noindent\textsc{Proof.}  Since cohomology groups, in particular the center, can not increase under deformation, it is sufficient to show that if $q$ is a primitive $N$th root of unity then $x^N, y^N$ are central.   Since $N_r = 0$ when $q$ is a primitive $N$th root of unity, (\ref{induction}) shows that in that case $x^N$ is central, and similarly for $y^N$.  $\Box$

\begin{theorem} The infinitesimal of the deformation from $W_{\text{qp}}$ to $W_q(\h)$ is $z_{\text{qp}}$.
\end{theorem}
\noindent\textsc{Proof.} Every element of $W_q(\h)$ can be written as a ``pseudopolynomial'', a linear combination over $k[\h]$ of element of the form $x^iy^j$.  Viewing $W_q(\h)$ as a deformation of $W_{\text{qp}}$, denote the multiplication in $W_q(\h)$ for the moment by `$\ast$'. To compute a product of the form $x^iy^m\ast x^ny^j$ one needs to know  $y^m\ast x^n$, but to know the infinitesimal only  the value modulo $\h^2$ is needed. By induction one finds
$$y^m\ast x^n \equiv r^mnx^my^n - r^{(m-1)(n-1)}m_rn_rr\h x^{m-1}y^{n-1} \pmod{\h^2}.$$ The bidegree of the right side is $(m-1,n-1)$ so the infinitesimal must be homogeneous of bidegree $(-1,-1)$. The only such non-trivial class is that of $z_{\text{qp}}$. $\Box$

\begin{theorem}\label{obstructed}   The derivations $x\partial_x$ and $y\partial_y$ of $W_{\text{qp}}$  are obstructed in the deformation of $W_{\text{qp}}$ to $W_q(\h)$ and can not be lifted to derivations of  $W_q(\h)$.   The derivation $x\partial_x-y\partial_y$ of $W_{\text{qp}}$ lifts to itself in $W_q(\h)$.  If a  derivation of  $W_{\text{qp}}$ which is not inner lifts to a derivation of $W_q(\h)$ then its lift is not inner.
\end{theorem}
\noindent\textsc{Proof.}  Since $z_{\text{qp}}$ is the infinitesimal of the deformation, the obstruction to lifting $x\partial_x$ is the class of $[x\partial_x, z_{\text{qp}}]$. Now $z_{\text{qp}}$ is itself a lift of $\partial_x\wedge\partial_y$, so the  class of $[x\partial_x, z_{\text{qp}}]$ is the class of a lift  of  $[x\partial_x,  \partial_x\wedge\partial_y]$. But the later is just $-\partial_x\wedge\partial_y$, which lifts to  $- z_{\text{qp}}$ and the class of the latter is not zero.  That $x\partial_x-y\partial_y$ lifts is evident since it preserves the defining equation. Finally, observe that no non-inner derivation of  $W_{\text{qp}}$, all of which are known from Theorem \ref{qp}, annihilates the entire center of $W_{\text{qp}}$, which is the same as that of $W_q(\h)$. Therefore no lift of a non-inner derivation can be inner by Theorem \ref{central}, $\Box$

\medskip
In the preceding theorem one can now set $\h = 1$.  When $q$ is not a root of unity it already gives the complete cohomology of $W_q$.

\begin{theorem} If $q$ is not an $N$th root of unity for any $N$ then $H^0(W_q) = k$, $H^1(W_q)$ is one-dimensional and generated by the class of $x\partial_x-y\partial_y$, and $H^2(W_q)$ is one-dimensional and generated by a class $w_q$ which is a lift of the class of $x\partial_x\wedge y\partial_y$.
\end{theorem}
\noindent\textsc{Proof.} We have seen that when $q$ is not a root of unity then $H^*(W_{\text{qp}})$ is spanned by the classes of $1, x\partial_x, y\partial_y, z _{\text{qp}}$ and $x\partial_x\wedge y\partial_y$; the \EP  characteristic is 1.  Since $z _{\text{qp}}$ lifts to a coboundary and the characteristic must remain equal to 1 after deformation it must be that $x\partial_x\wedge y\partial_y$ lifts but does not become a coboundary. $\Box$

\medskip
Since $H^3  = 0$ every 2-cocycle lifts; the essential statement about $w_q$ in the foregoing is that its class is not zero.   The class of $w_q$ is not zero in all cases, an essential fact that will be shown directly.  The preceding proof does not carry over to the case where $q$ is a root of unity because the cohomology of $W_{\text{qp}}$ is too large and the deformation has destroyed the grading, so the Hochschild complex can not be broken up into finite pieces.  It would be desirable to be able to prove the  result by a unifying argument on group actions similar to that of \S\ref{gradations}. (Replacing $q$ by $q^{\exp \hbar}$ would allow an action of the additive group but one must prove that this in fact gives a family of deformations.)  

Henceforth we assume that $q$ is a primitive $N$th root of unity. 

\begin{theorem}  The 2-cocycle $x\partial_x\wedge y\partial$ of $W_{\text{qp}}$ lifts to a 2-cocycle $w_q$ of $W_q$ and no multiple of  $w_q$ by an non-zero central element  is a coboundary.  The derivations $x^Ny\partial_y, y^Nx\partial_x, x^Nx\partial_x, y^Ny\partial_y$ all lift from $W_{\text{qp}}$ to $W_q(\h)$ and no non-zero element of the module they span over the center is a coboundary.
\end{theorem}
\noindent\textsc{Proof.}  Let  $c$ be any non-zero central element.
That $cx\partial_x\wedge y\partial$ lifts follows, as mentioned, from the fact that $H^3 = 0$. It follows, then, that $[cx\partial_x\wedge y\partial, x^N/N]= cx^Ny\partial_y$ also lifts (since $x^N$ lifts, in fact to itself). If the lift of $cx\partial_x\wedge y\partial$ were a coboundary then  $cx^Ny\partial_y$ would lift to a coboundary. But then, by Theorem \ref{central}, $cx^Ny\partial_y$ would have to annihilate every central element of $W_{\text{qp}}$ (since the all lift to central elements, namely themselves ), but this is clearly not the case. This also shows that $x^Ny\partial_y$ lifts, and since $x\partial_x-y\partial_y$ lifts, multiplying by $x^N$ shows that $x^Nx\partial_x$ also lifts.  The others follow similarly. To see that no non-zero element of the module these lifts span span can be a coboundary, observe that the same combination of $x^Ny\partial_y, y^Nx\partial_x, x^Nx\partial_x, y^Ny\partial_y$ (which is what was lifted) would again have to annihilate every element of the center, which is not the case. $\Box$

\medskip
Setting $\h = 1$ completely determines the cohomology of  $W_q$. Its center is $k[x^N, y^N]$ and $H^2(W_q)$ is a free module of rank 1 over the center generated by the class of $w_q$.  However, $H^1(W_q)$ is not quite a full free module since at the lowest level it contains only the class of $x\partial_x-y\partial_y$ but not lifts of $x\partial_x$ or $y\partial_y$. To generate $H^1(W_q)$ one needs the classes of $x\partial_x-y\partial_y$ and those of the lifts of $x^Nx\partial_x$ and $y^Ny\partial_y$. (Since the class of $x\partial_x-y\partial_y$ is present, those of $x^Ny\partial_y$ and $y^Nx\partial_x$ could be used instead.) In summary, we have the following.

\begin{theorem} When $q$ is a primitive $N$th root of unity then the center of $W_q$ is $k[x^N,y^N]$, $H^1(W_q)$ is the submodule of the free module over the center spanned by the classes of $x\partial_x-y\partial_y$ and the lifts of $x^Nx\partial_x, y^N\partial_y$, and $H^2(W_q)$ is the free module over the center spanned by the class of the lift $w_q$ of $x\partial_x\wedge y\partial_y$.$\Box$
\end{theorem}

Deformation theory shows that when $q$ is a root of unity then certain classes must lift  but does not exhibit the lifts explicitly. That they exist depends on the fact that $q$ is a root of unity and therefore encodes some information about roots of unity. It might be no more than that the sum of the roots vanishes (a fact used here explicitly), but it could be interesting to know the actual lifts.

\end{document}